\documentclass[10pt]{amsart}
\usepackage{latexsym}   
\usepackage{amssymb}    
\usepackage{amsmath}    
\usepackage{amsthm}     
\usepackage[all]{xy}
\usepackage{hyperref}
\newcommand{\pa}{\partial}\newcommand{\al}{\alpha}

\newcommand{\om}{\omega}

\newcommand{\ti}{\tilde}

\renewcommand{\thefootnote}

\newtheorem{theorem}{Theorem}[section]

\theoremstyle{definition}
\newtheorem{definition}[theorem]{Definition}

\theoremstyle{remark}

\numberwithin{equation}{section}

\title[On isometric correspondence of leaves]
{On isometric correspondence of leaves}

\author[  Ion I. Dinc\u{a}]{Ion I. Dinc\u{a}}
\address{Department of Applied Mathematics,
Faculty of Applied Sciences, National University of Science and Technology Politehnica Bucharest
313 Spl.Independentei 060042 Bucharest, Romania}
 \email{idinca@upb.ro}

\begin{document}


\begin{abstract}
We prove that for a generic $3$-dimensional integrable rolling distribution of contact elements (excluding developable seed and isotropic developable leaves) isometric correspondence of leaves of a general nature (independent of the shape of the seed) requires the B\"{a}cklund transformation.
\end{abstract}

\maketitle

\tableofcontents \pagenumbering{arabic}

\section{Introduction}

We shall consider the complexification
$$(\mathbb{C}^3,\langle\cdot,\cdot\rangle),\ \langle x,y\rangle :=x^Ty,\
|x|^2:=x^Tx,\ x,y\in\mathbb{C}^3$$
of the real $3$-dimensional Euclidean space; in this setting surfaces are $2$-dimensional objects of
$\mathbb{C}^3$ depending on two real or complex parameters.

{\it Isotropic} (null) vectors are those vectors of length $0$; since most vectors are not isotropic
we call a vector simply vector and we shall emphasize isotropic for isotropic vectors. The same
denomination will apply in other settings: for example we call quadric a non-degenerate quadric (a
quadric projectively equivalent to the complex unit sphere).

Consider Lie's viewpoint: one can replace a surface $x\subset\mathbb{C}^3$ with a $2$-dimensional
distribution of {\it contact elements} (pairs of points and planes passing through those points; the
classical geometers call them {\it facets}): the collection of its tangent planes (with the points
of tangency highlighted); thus a contact element is the infinitesimal version of a surface (the
integral element $(x,dx)|_{\mathrm{pt}}$ of the surface). Conversely, a $2$-dimensional distribution
of contact elements is not always the collection of the tangent planes of a surface (with the points
of tangency highlighted), but the condition that a $2$-dimensional distribution of contact elements
is integrable (that is it is the collection of the tangent planes of a leaf (sub-manifold)) does not
distinguish between the cases when this sub-manifold is a surface, curve or point, thus allowing the
collapsing of the leaf.

A $3$-dimensional distribution of contact elements is integrable if it is the collection of the
tangent planes of an $1$-dimensional family of leaves.

Two {\it rollable} (applicable or isometric) surfaces can be {\it rolled} (applied) one onto the
other such that at any instant they meet tangentially and with same differential at the tangency
point.

\begin{definition}
The rolling of two isometric surfaces $x_0,\ x\subset\mathbb{C}^3$ (that is $|dx_0|^2=|dx|^2$) is
the surface, curve or point $(R,t)\subset\mathbf{O}_3(\mathbb{C})\ltimes\mathbb{C}^3$ such that
$(x,dx)=(R,t)(x_0,dx_0):=(Rx_0+t,Rdx_0)$.
\end{definition}

The rolling introduces the flat connection form (it encodes the difference of the second fundamental
forms of $x_0,\ x$ and it being flat encodes the difference of the Gau\ss-Codazzi-Mainardi-Peterson equations of $x_0,\ x$).

\begin{definition}
Consider an integrable $3$-dimensional distribution of contact elements $\mathcal{F}=(p,m)$ centered
at $p=p(u,v,w)$, with normal fields $m=m(u,v,w)$ and distributed along the surface $x_0=x_0(u,v)$.
If we roll $x_0$ on an isometric surface $x$ (that is $(x,dx)=(R,t)(x_0,dx_0):=(Rx_0+t,Rdx_0)$), then
the rolled distribution of contact elements is $(Rp+t,Rm)$ and is distributed along $x$; if it
remains integrable for any rolling, then the distribution is called {\it $3$-dimensional integrable rolling
distribution of contact elements} with seed $x$ and leaf $Rp+t$.
\end{definition}

Bianchi considered the most general form of a B\"{a}cklund transformation as the focal surfaces (one
transform of the other) of a {\it Weingarten} congruence (congruence upon whose two focal surfaces
the asymptotic directions correspond; equivalently the second fundamental forms are proportional).
Note that although the correspondence provided by the Weingarten congruence does not give the
applicability (isometric) correspondence, the B\"{a}cklund transformation is the tool best suited to
attack the isometric deformation problem via geometric transformation, since it provides correspondence of the
characteristics of the isometric deformation problem (according to Darboux these are the asymptotic
directions), it is directly linked to the infinitesimal isometric deformation problem (Darboux
proved that infinitesimal isometric deformations generate Weingarten congruences and Guichard proved
the converse: there is an infinitesimal isometric deformation of a focal surface of a Weingarten
congruence in the direction normal to the other focal surface; see Darboux (\cite{D1},\S\ 883-\S\
924)) and it admits a version of the Bianchi Permutability Theorem for its second iteration.

With $\ti d\cdot:=\pa_u\cdot du+\pa_v\cdot dv+\pa_w\cdot dw=d\cdot+\pa_w\cdot dw, V:=p-x_0$ if the $3$-dimensional distribution
of contact elements is integrable and the rolled distribution remains integrable if we roll $x_0$ on
an isometric surface $x,\ (x,dx)=(R,t)(x_0,dx_0)$ (that is we replace $x_0,\ V,\ m$ with $x,\ RV,\
Rm$), then along the leaves we have

\begin{eqnarray}\label{eq:int}
0=(Rm)^T\ti d(RV+x)=m^T(\om\times V+d(V+x_0)+\pa_wVdw),\nonumber\\ \om:=N_0\times R^{-1}dRN_0,\ N_0:=\frac{\pa_ux_0\times\pa_vx_0}{|\pa_ux_0\times\pa_vx_0|}.
\end{eqnarray}
Since we shall not need the integrability condition of this general integrable rolling distribution
of contact elements, we shall not derive it.
\begin{definition}
A $3$-dimensional integrable rolling distribution of contact elements is called generic if $m^T\pa_wV(m\times V)\times N_0\neq 0$.
\end{definition}

In \cite{D3} we proved that for a generic $3$-dimensional integrable rolling distribution of contact elements (excluding developable seed and isotropic developable leaves) and with the symmetry of the tangency configuration (contact elements are centered on tangent planes of the surface $x_0$ and further pass through the origin of the tangent planes) the seed and any leaf are the focal surfaces of a Weingarten congruence (and thus we get B\"{a}cklund transformation according to Bianchi's definition) and for a generic $3$-dimensional integrable rolling distribution of contact elements (excluding developable seed and isotropic developable leaves) isometric correspondence of leaves of a general nature (independent of the shape of the seed) requires the tangency configuration (contact elements are centered on tangent planes of the surface $x_0$).
Further by applying, if necessary, a change of variables $w=w(\ti w, u,v)$ we have
$N_0^T[d(V+x_0)\times\wedge d(V+x_0)]\neq 0$ and we get the condition
\begin{eqnarray}\label{eq:cons}
d(V+x_0)^T(I_3-\frac{mN_0^T}{m^TN_0})\odot[d(\frac{N_0m^T}{m^TN_0})\pa_wV-
\pa_w(\frac{N_0m^T}{m^TN_0})d(V+x_0)+\nonumber\\
2N_0^T[\pa_wV\times d(V+x_0)]\frac{d(\frac{N_0m^T}{m^TN_0})\wedge
d(V+x_0)}{N_0^T[d(V+x_0)\times\wedge d(V+x_0)]}]=0,\ du\wedge dv\wedge dw\neq 0.
\end{eqnarray}

According to Bianchi \cite{B3} (referring to $3$-dimensional integrable rolling distributions of contact elements) {\it "But, in view of the eventual applications to problems of deformations, it is opportune to limit the problem much more, and to suppose that every facette $f$ and each of its associated facettes $f'$ has the center of one in the plane of the other"}.

We have now the main {\bf Theorem} of this paper:

\begin{theorem}\label{th:th1}
For a generic $3$-dimensional integrable rolling distribution of contact elements (excluding developable seed and isotropic developable leaves) isometric correspondence of leaves of a general nature (independent of the shape of the seed) requires the B\"{a}cklund transformation.
\end{theorem}

The remaining part of the paper is organized as follows: in Section 2 we recall the rolling problem for surfaces and in Section 3 we provide the proof of {\bf Theorem} \ref{th:th1}.

\section{The rolling problem for surfaces}

Let $(u,v)\in D$ with $D$ domain of $\mathbb{R}^2$ or $\mathbb{C}^2$ and $x:D\mapsto\mathbb{C}^3$ be
a surface.

For $\om_1,\om_2\ \mathbb{C}^3$-valued $1$-forms on $D$ and $a,b\in\mathbb{C}^3$ we have
\begin{eqnarray}\label{eq:fund}
a^T\om_1\wedge b^T\om_2=((a\times b)\times\om_1+b^T\om_1a)^T\wedge\om_2=
(a\times b)^T(\om_1\times\wedge\om_2)+b^T\om_1\wedge a^T\om_2;\nonumber\\
\mathrm{in\ particular}\ a^T\om\wedge b^T\om=\frac{1}{2}(a\times b)^T(\om\times\wedge\om).
\end{eqnarray}

Since both $\times$ and $\wedge$ are skew-symmetric, we have $2\om_1\times\wedge\om_2=
\om_1\times\om_2+\om_2\times\om_1=2\om_2\times\wedge\om_1$.

Consider the scalar product $\langle\cdot,\cdot\rangle$ on $\mathbf{M}_3(\mathbb{C}):\
\langle X,Y\rangle :=\frac{1}{2}\mathrm{tr}(X^TY)$. We have the isometry

$$\al:\mathbb{C}^3\mapsto\mathbf{o}_3(\mathbb{C}),\
\al\left(\begin{bmatrix}x^1\\x^2\\x^3\end{bmatrix}\right)
=\begin{bmatrix}0&-x^3&x^2\\x^3&0&-x^1\\-x^2&x^1&0\end{bmatrix},\
x^Ty=\langle\al(x),\al(y)\rangle=
\frac{1}{2}\mathrm{tr}(\al(x)^T\al(y)),$$
$$\al(x\times y)=[\al(x),\al(y)]=\al(\al(x)y)=yx^T-xy^T,\ \al(Rx)=R\al(x)R^{-1},\
x,y\in\mathbb{C}^3,\ R\in\mathbf{O}_3(\mathbb{C}).$$

Let $x\subset\mathbb{C}^3$ be a surface applicable (isometric) to a surface
$x_0\subset\mathbb{C}^3$:

\begin{eqnarray}\label{eq:roll}
(x,dx)=(R,t)(x_0,dx_0):=(Rx_0+t,Rdx_0),
\end{eqnarray}
where $(R,t)$ is a sub-manifold in $\mathbf{O}_3(\mathbb{C})\ltimes\mathbb{C}^3$ (in general
surface, but it is a curve if $x_0,\ x$ are ruled and the rulings correspond under isometry or a
point if $x_0,\ x$ differ by a rigid motion). The sub-manifold $R$ gives the rolling of $x_0$ on $x$,
that is if we rigidly roll $x_0$ on $x$ such that points corresponding under the isometry will have
the same differentials, $R$ will dictate the rotation of $x_0$; the translation $t$ will satisfy
$dt=-dRx_0$.

For $(u,v)$ parametrization on $x_0,\ x$ and outside the locus of isotropic (degenerate) induced
metric of $x_0,\ x$ we have $N_0:=\frac{\pa_ux_0\times\pa_vx_0}{|\pa_ux_0\times\pa_vx_0|},\
N:=\frac{\pa_ux\times\pa_vx}{|\pa_ux\times\pa_vx|}$ respectively positively oriented unit normal
fields of $x_0,\ x$ and $R$ is determined by $R=[\pa_ux\ \ \pa_vx\ \ N][\pa_ux_0\ \ \pa_vx_0\ \
\det(R)N_0]^{-1}$; we take $R$ with $\det(R)=1$; thus the rotation of the rolling with the other
face of $x_0$ (or on the other face of $x$) is $R':=R(I_3-2N_0N_0^T)=(I_3-2NN^T)R,\ \det(R')=-1$.

Therefore $\mathbf{O}_3(\mathbb{C})\ltimes\mathbb{C}^3$ acts on
$2$-dimensional integrable distributions of contact elements $(x_0,dx_0)$ in
$T^*(\mathbb{C}^3)$ as: $(R,t)(x_0,dx_0)=(Rx_0+t,Rdx_0)$; a
rolling is a sub-manifold
$(R,t)\subset\mathbf{O}_3(\mathbb{C})\ltimes\mathbb{C}^3$ such
that $(R,t)(x_0,dx_0)$ is still integrable.

We have:
\begin{eqnarray}\label{eq:secoi}
R^{-1}dRN_0=R^{-1}dN-dN_0.
\end{eqnarray}
Applying the
compatibility condition $d$ to (\ref{eq:roll}) we get:
\begin{eqnarray}\label{eq:comp}
R^{-1}dR\wedge dx_0=0,\ dRR^{-1}\wedge dx=0.
\end{eqnarray}

Since $R^{-1}dR$ is skew-symmetric and using (\ref{eq:comp}) we
have
\begin{eqnarray}\label{eq:dx0}
dx_0^TR^{-1}dRdx_0=0.
\end{eqnarray}
From (\ref{eq:dx0}) for $a\in\mathbb{C}^3$ we get
$R^{-1}dRa=R^{-1}dR(a^{\bot}+a^{\top})=a^TN_0R^{-1}dRN_0-a^TR^{-1}dRN_0N_0=\om\times
a,\ \om:=N_0\times R^{-1}dRN_0=^{(\ref{eq:secoi})}(\det
R)R^{-1}(N\times dN)-N_0\times dN_0=R^{-1}(N\times dN)-N_0\times
dN_0$. Thus $R^{-1}dR=\al(\om)$ and $\om$ is flat connection form
in $T^*x_0$:
\begin{eqnarray}\label{eq:om}
d\om+\frac{1}{2}\om\times\wedge\om=0,\ \om\times\wedge
dx_0=0,\ (\om)^\perp=0.
\end{eqnarray}
With $s:=N_0^T(\om\times
dx_0)=s_{11}du^2+s_{12}dudv+s_{21}dvdu+s_{22}dv^2$ the difference
of the second fundamental forms of $x,\ x_0$ we have
\begin{eqnarray}\label{eq:omjk}
\om=\frac{s_{12}\pa_ux_0-s_{11}\pa_vx_0}{|\pa_ux_0\times\pa_vx_0|}du+
\frac{s_{22}\pa_ux_0-s_{21}\pa_vx_0}{|\pa_ux_0\times\pa_vx_0|}dv;
\end{eqnarray}
($\om\times\wedge dx_0=0$ is equivalent to $s_{12}=s_{21}$;
$(d\om)^\perp+\frac{1}{2}\om\times\wedge\om=0,\
(d\om)^\top=0$ respectively encode the difference of the
Gau\ss\ -Codazzi-Mainardi-Peterson equations of $x_0$ and $x$).

Using $\frac{1}{2}dN_0\times\wedge
dN_0=K|\pa_ux_0\times\pa_vx_0|N_0du\wedge dv,\ K$ being the
Gau\ss\ curvature we get $dN_0\times\wedge
dN_0=R^{-1}(dN\times\wedge dN)=^{(\ref{eq:secoi})}(\om\times
N_0+dN_0)\times\wedge(\om\times N_0+dN_0)=dN_0\times\wedge
dN_0+2(\om\times N_0)\times\wedge dN_0+\om\times\wedge\om$; thus
\begin{eqnarray}\label{eq:omom}
\frac{1}{2}\om\times\wedge\om=dN_0^T\wedge\om N_0.
\end{eqnarray}
Note also
\begin{eqnarray}\label{eq:om'}
\om'=N_0\times {R'}^{-1}dR'N_0=-\om-2N_0\times dN_0
\end{eqnarray}
and
\begin{eqnarray}\label{eq:aom}
\ \ \ \ a^T\wedge\om=0,\ \forall\om\ \mathrm{satisfying\
(\ref{eq:om})\ for}\ a\ 1-\mathrm{form}\Rightarrow\ a^T\odot
dx_0:=\frac{a^Tdx_0+dx_0^Ta}{2}=0.
\end{eqnarray}
Note that the converse $a^T\odot dx_0^0=0,\ a\ 1-$form
$\Rightarrow a^T\wedge\om=0,\forall\om$ satisfying (\ref{eq:om})
is also true.

\section {Proof of Theorem \ref{th:th1}}

\subsection{The case $m=V\times N_0+\mathbf{m}N_0+\mathbf{n}V$}\noindent

\noindent

We have $0=(Rm)^T\ti d(RV+x)=m^T(\om\times V+d(V+x_0)+\pa_wVdw)$, or, assuming $m^T\pa_wV\neq 0$,

$$dw=\frac{N_0^T[V\times d(V+x_0)]+\mathbf{m}V^T(\om\times N_0+dN_0)-\mathbf{n}V^Td(V+x_0)}
{N_0^T(\pa_wV\times V)+\mathbf{n}V^T\pa_wV}.$$

Applying the compatibility condition $d$ to this equation and using the equation itself we get

$$0=-\frac{\pa_w[N_0^T[V\times d(V+x_0)]+\mathbf{m}V^T(\om\times N_0+dN_0)-\mathbf{n}V^Td(V+x_0)]}
{N_0^T(\pa_wV\times V)+\mathbf{n}V^T\pa_wV}\wedge$$
$$\frac{N_0^T[V\times d(V+x_0)]+\mathbf{m}V^T(\om\times N_0+dN_0)-\mathbf{n}V^Td(V+x_0)}
{N_0^T(\pa_wV\times V)+\mathbf{n}V^T\pa_wV}+$$
$$d\frac{N_0^T[V\times d(V+x_0)]}{N_0^T(\pa_wV\times V)+\mathbf{n}V^T\pa_wV}+
d\frac{\mathbf{m}V^T}{N_0^T(\pa_wV\times V)+\mathbf{n}V^T\pa_wV}\wedge(\om\times N_0+dN_0)-$$
$$d\frac{\mathbf{n}V^T}{N_0^T(\pa_wV\times V)+\mathbf{n}V^T\pa_wV}\wedge d(V+x_0)$$
$$\stackrel{(\ref{eq:fund})}{=}-\frac{\pa_w[N_0^T[V\times d(V+x_0)]-\mathbf{n}V^Td(V+x_0)]}
{N_0^T(\pa_wV\times V)+\mathbf{n}V^T\pa_wV}\wedge\frac{N_0^T[V\times d(V+x_0)]-\mathbf{n}V^Td(V+x_0)}
{N_0^T(\pa_wV\times V)+\mathbf{n}V^T\pa_wV}+$$
$$d\frac{N_0^T[V\times d(V+x_0)]}{N_0^T(\pa_wV\times V)+\mathbf{n}V^T\pa_wV}-
d\frac{\mathbf{n}V^T}{N_0^T(\pa_wV\times V)+\mathbf{n}V^T\pa_wV}\wedge d(V+x_0)-$$
$$\frac{\mathbf{m}^2N_0^T(\pa_wV\times V)N_0^T(dN_0\times\wedge dN_0)}{2[N_0^T(\pa_wV\times V)+\mathbf{n}V^T\pa_wV]^2}+[\frac{\pa_w\mathbf{m}[N_0^T[V\times d(V+x_0)] -\mathbf{n}V^Td(V+x_0)]}{[N_0^T(\pa_wV\times V)+\mathbf{n}V^T\pa_wV]^2}-$$
$$\frac{\mathbf{m}\pa_w[N_0^T[V\times d(V+x_0)]-\mathbf{n}V^Td(V+x_0)]}{[N_0^T(\pa_wV\times V)+\mathbf{n}V^T\pa_wV]^2}+d\frac{\mathbf{m}}{N_0^T(\pa_wV\times V)+\mathbf{n}V^T\pa_wV}-$$
$$\frac{\mathbf{m}(\pa_wV\times N_0)^Td(V+x_0)}{[N_0^T(\pa_wV\times V)+\mathbf{n}V^T\pa_wV]^2}]\wedge V^T(\om\times N_0+dN_0)+$$
$$\frac{\mathbf{m}\mathbf{n}N_0^T[\pa_wV\times d(V+x_0)]}{[N_0^T(\pa_wV\times V)+\mathbf{n}V^T\pa_wV]^2}
\wedge(N_0\times V)^T(\om\times N_0+dN_0)$$
$$\stackrel{(\ref{eq:fund})}{=}\frac{1}{[N_0^T(\pa_wV\times V)+\mathbf{n}V^T\pa_wV]^2}[-|V|^2d\mathbf{n}\wedge
(N_0\times\pa_wV)^Td(V+x_0)-\frac{1}{2}[N_0^T(\pa_wV\times V)(1+\mathbf{n}^2)$$
$$-\pa_w\mathbf{n}|V|^2]N_0^T[d(V+x_0)\times\wedge d(V+x_0)]-\frac{1}{2}|V|^2[N_0^T(\pa_wV\times V)+\mathbf{n}V^T\pa_wV]N_0^T(dN_0\times\wedge dN_0)$$
$$+[N_0^T(\pa_wV\times V)+\mathbf{n}V^T\pa_wV]N_0^T[dV\times\wedge d(V+x_0)]-(N_0\times\pa_wV+\mathbf{n}\pa_wV)^TdV\wedge(N_0\times V)^Td(V+x_0)$$
$$+\mathbf{n}[(dV\times V)\times(N_0\times\pa_wV+\mathbf{n}\pa_wV)]^T\wedge d(V+x_0)-\frac{\mathbf{m}^2}{2}N_0^T(\pa_wV\times V)N_0^T(dN_0\times\wedge dN_0)]$$
$$+\frac{1}{N_0^T(\pa_wV\times V)+\mathbf{n}V^T\pa_wV}[d\mathbf{m}+\frac{\pa_w\mathbf{m}[N_0^T[V\times d(V+x_0)] -\mathbf{n}V^Td(V+x_0)]}{N_0^T(\pa_wV\times V)+\mathbf{n}V^T\pa_wV}$$
$$+\frac{\mathbf{m}[\pa_w\mathbf{n}V^Td(V+x_0)+\mathbf{n}\pa_wV^Td(V+x_0)]}{N_0^T(\pa_wV\times V)+\mathbf{n}V^T\pa_wV}-\frac{\mathbf{m}[N_0^T(\pa_wV\times dV)+d\mathbf{n}V^T\pa_wV+\mathbf{n}dV^T\pa_wV]}{N_0^T(\pa_wV\times V)+\mathbf{n}V^T\pa_wV}]\wedge$$
$$V^T(\om\times N_0+dN_0)+\frac{\mathbf{m}\mathbf{n}N_0^T[\pa_wV\times d(V+x_0)]}{[N_0^T(\pa_wV\times V)+\mathbf{n}V^T\pa_wV]^2}
\wedge(N_0\times V)^T(\om\times N_0+dN_0).$$

This can be written for short

\begin{eqnarray}\label{eq:Vom}
\ \ \ \ \ A+B\wedge V^T(\om\times N_0+dN_0)+C\wedge(N_0\times V)^T(\om\times N_0+dN_0)=0,\ \forall\ \om\ \mathrm{satisfying}\ (\ref{eq:om}),
\end{eqnarray}
where $A$ is a scalar $2$-form not depending on $\om$ and $B,C$ are scalar $1$-forms not depending on $\om$.

In order for $w$ to be determined by $\om$ from (\ref{eq:int}), we need $(m\times V)\times N_0\neq 0$ and $w$ cannot be linked to $\om$ by any other relation, either functional (as (\ref{eq:Vom}) a-priori is) or differential, thus
in (\ref{eq:Vom}) $\om$ cancels independently of $w$ and outside $w$ we can replace $\om$ with any other solution of (\ref{eq:om}).

Replacing $\om$ respectively with $0,-2N_0\times dN_0,-\om-2N_0\times dN_0$ we get $A=0$ and
\begin{eqnarray}\label{eq:BC}
(BV+CN_0\times V)^T\wedge dN_0=0,\ (BN_0\times V-CV)^T\wedge\om=0,\ \forall\ \om\ \mathrm{satisfying}\ (\ref{eq:om}).
\end{eqnarray}
From (\ref{eq:aom}) we get from the second relation of (\ref{eq:BC}) $(BN_0\times V-CV)^T\odot dx_0=0$; with $B=B_udu+B_vdv,\ C=C_udu+C_vdv$ this becomes $B_u(N_0\times V)^T\pa_ux_0-C_uV^T\pa_ux_0=B_v(N_0\times V)^T\pa_vx_0-C_vV^T\pa_vx_0=B_u(N_0\times V)^T\pa_vx_0+B_v(N_0\times V)^T\pa_ux_0-C_uV^T\pa_vx_0-C_vV^T\pa_ux_0=0$, so $B_u=C_u\frac{V^T\pa_ux_0}{(N_0\times V)^T\pa_ux_0},\
B_v=C_v\frac{V^T\pa_vx_0}{(N_0\times V)^T\pa_vx_0}$ and $(N_0\times V)^T(C_u\pa_vx_0-C_v\pa_ux_0)=0$.
With $V=V_1\pa_ux_0+V_2\pa_vx_0$ we get $C_v=-aV_1,\ C_u=aV_2,\ a\subset\mathbb{C}$.

In the first relation of (\ref{eq:BC}) we can replace $x_0$ with any isometric surface $x$; thus
$BV+CN_0\times V=bdx_0,\ b\subset\mathbb{C}$ and $C=\frac{bdx_0^T(N_0\times V)}{|N_0\times V|^2}$, which is equivalent to the previous determination of $C$.

We thus get $$\frac{\mathbf{m}\mathbf{n}N_0^T[\pa_wV\times d(V+x_0)]}{[N_0^T(\pa_wV\times V)+\mathbf{n}V^T\pa_wV]^2}=\frac{bdx_0^T(N_0\times V)}{|N_0\times V|^2}$$ and since $\mathbf{m}\neq 0$ we get

\begin{eqnarray}\label{eq:n}
\mathbf{n}N_0^T[\pa_wV\times d(V+x_0)]\wedge dx_0^T(N_0\times V)=0
\end{eqnarray}
and

\begin{eqnarray}\label{eq:b}
b=\frac{\mathbf{m}\mathbf{n}|N_0\times V|^2}{[N_0^T(\pa_wV\times V)+\mathbf{n}V^T\pa_wV]^2}\frac{N_0^T[\pa_wV\times\pa_u(V+x_0)]}{\pa_ux_0^T(N_0\times V)}\nonumber\\
=\frac{\mathbf{m}\mathbf{n}|N_0\times V|^2}{[N_0^T(\pa_wV\times V)+\mathbf{n}V^T\pa_wV]^2}\frac{N_0^T[\pa_wV\times\pa_v(V+x_0)]}{\pa_vx_0^T(N_0\times V)}.
\end{eqnarray}

Now equation (\ref{eq:cons}) becomes:

$$\frac{d(V+x_0)^T\odot dN_0}{\mathbf{m}}[N_0^T(\pa_wV\times V)+\mathbf{n}V^T\pa_wV]-\frac{1}{\mathbf{m}}d(V+x_0)^T(V\times N_0+\mathbf{n}V)\odot[-\frac{d\mathbf{m}}{\mathbf{m}^2}[N_0^T(\pa_wV\times V)$$
$$+\mathbf{n}V^T\pa_wV]+\frac{1}{\mathbf{m}}(dV\times N_0+\mathbf{m}dN_0+d\mathbf{n}V+\mathbf{n}dV)^T\pa_wV
+\frac{\pa_w\mathbf{m}}{\mathbf{m}^2}(V\times N_0+\mathbf{n}V)^Td(V+x_0)$$
$$-\frac{1}{\mathbf{m}}(\pa_wV\times N_0+\pa_w\mathbf{n}V+\mathbf{n}\pa_wV)^Td(V+x_0)]+\frac{N_0^T[\pa_wV\times d(V+x_0)]}{N_0^T[\pa_u(V+x_0)\times\pa_v(V+x_0)]}\odot$$
$$[\frac{d(V+x_0)^T\pa_uN_0}{\mathbf{m}}(V\times N_0+\mathbf{m}N_0+\mathbf{n}V)^T\pa_v(V+x_0)-\frac{d(V+x_0)^T\pa_vN_0}{\mathbf{m}}(V\times N_0+\mathbf{m}N_0$$$$+\mathbf{n}V)^T\pa_u(V+x_0)-\frac{1}{\mathbf{m}}d(V+x_0)^T(V\times N_0+\mathbf{n}V)[-\frac{\pa_u\mathbf{m}}{\mathbf{m}^2}(V\times N_0+\mathbf{n}V)^T\pa_v(V+x_0)$$
$$+\frac{\pa_v\mathbf{m}}{\mathbf{m}^2}(V\times N_0+\mathbf{n}V)^T\pa_u(V+x_0)+\frac{1}{\mathbf{m}}(\pa_uV\times N_0+V\times\pa_uN_0+\mathbf{m}\pa_uN_0+\pa_u\mathbf{n}V+\mathbf{n}\pa_uV)^T\pa_v(V+x_0)$$
$$-\frac{1}{\mathbf{m}}(\pa_vV\times N_0+V\times\pa_vN_0+\mathbf{m}\pa_vN_0+\pa_v\mathbf{n}V+\mathbf{n}\pa_vV)^T\pa_u(V+x_0)]]=0.$$

Taking $d\mathbf{m}$ from $B=\frac{bdx_0^TV}{|V|^2}$ and using $A=0$ the derivatives of $\mathbf{m}$ and $\mathbf{n}$ in the above equation disappear and we get only:

$$\mathbf{n}\{\mathbf{n}^2[\frac{N_0^T[\pa_wV\times d(V+x_0)]}{N_0^T[\pa_u(V+x_0)\times\pa_v(V+x_0)]}\odot\frac{d(V+x_0)^TV}{\mathbf{m}^2}
[N_0^T(\pa_wV\times V)N_0^T[\pa_u(V+x_0)\times\pa_v(V+x_0)]-$$$$[(\pa_uV\times V)\times\pa_wV]^T\pa_v(V+x_0)+
[(\pa_vV\times V)\times\pa_wV]^T\pa_u(V+x_0)]+$$$$\frac{1}{\mathbf{m}^2}d(V+x_0)^TV\odot dx_0^TV\frac{N_0^T[\pa_wV\times\pa_u(V+x_0)]}{\pa_ux_0^T(N_0\times V)}V^T\pa_wV+$$$$
\frac{N_0^T[\pa_wV\times d(V+x_0)]}{N_0^T[\pa_u(V+x_0)\times\pa_v(V+x_0)]}\odot\frac{d(V+x_0)^TV}{\mathbf{m}^2}
[[\pa_ux_0^TV\frac{N_0^T[\pa_wV\times\pa_u(V+x_0)]}{\pa_ux_0^T(N_0\times V)}-\pa_wV^T\pa_u(V+x_0)+$$$$\pa_uV^T\pa_wV]V^T\pa_v(V+x_0)-
[\pa_vx_0^TV\frac{N_0^T[\pa_wV\times\pa_u(V+x_0)]}{\pa_ux_0^T(N_0\times V)}-\pa_wV^T\pa_v(V+x_0)+$$$$\pa_vV^T\pa_wV]V^T\pa_u(V+x_0)]]+\mathbf{n}[\frac{N_0^T[\pa_wV\times d(V+x_0)]}{N_0^T[\pa_u(V+x_0)\times\pa_v(V+x_0)]}\odot\frac{d(V+x_0)^TV}{\mathbf{m}^2}[
|V|^2V^T\pa_wVN_0^T(\pa_uN_0\times\pa_vN_0)$$$$-V^T\pa_wV[N_0^T[\pa_uV\times\pa_v(V+x_0)]
-N_0^T[\pa_vV\times\pa_u(V+x_0)]]+\pa_wV^T\pa_uV(N_0\times V)^T\pa_v(V+x_0)-$$$$\pa_wV^T\pa_vV(N_0\times V)^T\pa_u(V+x_0)-[(\pa_uV\times V)\times(N_0\times\pa_wV)]^T\pa_v(V+x_0)+[(\pa_vV\times V)\times(N_0\times\pa_wV)]^T\pa_u(V+x_0)]+$$$$\frac{N_0^T[\pa_wV\times d(V+x_0)]}{N_0^T[\pa_u(V+x_0)\times\pa_v(V+x_0)]}\odot\frac{d(V+x_0)^T(V\times N_0)}{\mathbf{m}^2}[N_0^T(\pa_wV\times V)N_0^T[\pa_u(V+x_0)\times\pa_v(V+x_0)]-$$$$[(\pa_uV\times V)\times\pa_wV]^T\pa_v(V+x_0)+
[(\pa_vV\times V)\times\pa_wV]^T\pa_u(V+x_0)]+$$$$\frac{d(V+x_0)^T\odot dN_0}{\mathbf{m}}(\pa_wV^TV)^2
+\frac{1}{\mathbf{m}^2}d(V+x_0)^TV\odot dx_0^TV\frac{N_0^T[\pa_wV\times\pa_u(V+x_0)]}{\pa_ux_0^T(N_0\times V)}
N_0^T(\pa_wV\times V)+$$$$\frac{1}{\mathbf{m}^2}d(V+x_0)^T(V\times N_0)\odot dx_0^TV\frac{N_0^T[\pa_wV\times\pa_u(V+x_0)]}{\pa_ux_0^T(N_0\times V)}V^T\pa_wV-\frac{1}{\mathbf{m}}d(V+x_0)^TV\odot dN_0^T\pa_wVV^T\pa_wV+$$$$\frac{1}{\mathbf{m}^2}d(V+x_0)^TV\odot (\pa_wV\times N_0)^Td(V+x_0)V^T\pa_wV+\frac{N_0^T[\pa_wV\times d(V+x_0)]}{N_0^T[\pa_u(V+x_0)\times\pa_v(V+x_0)]}\odot[
$$$$\frac{d(V+x_0)^T\pa_uN_0}{\mathbf{m}}V^T\pa_v(V+x_0)V^T\pa_wV
-\frac{d(V+x_0)^T\pa_vN_0}{\mathbf{m}}V^T\pa_u(V+x_0)V^T\pa_wV-$$$$
\frac{1}{\mathbf{m}}d(V+x_0)^T(V\times N_0)[(-\frac{1}{\mathbf{m}}\pa_ux_0^TV\frac{N_0^T[\pa_wV\times\pa_u(V+x_0)]}{\pa_ux_0^T(N_0\times V)}+\frac{1}{\mathbf{m}}\pa_wV^T\pa_u(V+x_0)-$$$$\frac{1}{\mathbf{m}}\pa_uV^T\pa_wV)V^T\pa_v(V+x_0)
+(\frac{1}{\mathbf{m}}\pa_vx_0^TV\frac{N_0^T[\pa_wV\times\pa_u(V+x_0)]}{\pa_ux_0^T(N_0\times V)}-\frac{1}{\mathbf{m}}\pa_wV^T\pa_v(V+x_0)+$$$$\frac{1}{\mathbf{m}}\pa_vV^T\pa_wV)V^T\pa_u(V+x_0)
+\frac{1}{\mathbf{m}}\pa_uV^T\pa_v(V+x_0)V^T\pa_wV-\frac{1}{\mathbf{m}}\pa_vV^T\pa_u(V+x_0)V^T\pa_wV]-$$$$
\frac{1}{\mathbf{m}}d(V+x_0)^TV[(-\frac{1}{\mathbf{m}}\pa_ux_0^TV\frac{N_0^T[\pa_wV\times
\pa_u(V+x_0)]}{\pa_ux_0^T(N_0\times V)}+\frac{1}{\mathbf{m}}\pa_wV^T\pa_u(V+x_0)-$$$$\frac{1}{\mathbf{m}}\pa_uV^T\pa_wV)(V\times N_0)^T\pa_v(V+x_0)+(\frac{1}{\mathbf{m}}\pa_vx_0^TV\frac{N_0^T[\pa_wV\times\pa_u(V+x_0)]}{\pa_ux_0^T(N_0\times V)}-\frac{1}{\mathbf{m}}\pa_wV^T\pa_v(V+x_0)+$$$$\frac{1}{\mathbf{m}}\pa_vV^T\pa_wV)(V\times N_0)^T\pa_u(V+x_0)-\frac{1}{\mathbf{m}}N_0^T(\pa_wV\times\pa_uV)V^T\pa_v(V+x_0)+$$$$
\frac{1}{\mathbf{m}}N_0^T(\pa_wV\times\pa_vV)V^T\pa_u(V+x_0)+
\frac{1}{\mathbf{m}}\pa_uV^T\pa_v(V+x_0)N_0^T(\pa_wV\times V)-$$$$
\frac{1}{\mathbf{m}}\pa_vV^T\pa_u(V+x_0)N_0^T(\pa_wV\times V)+\frac{1}{\mathbf{m}}(\pa_uV\times N_0+\mathbf{m}\pa_uN_0)^T\pa_v(V+x_0)V^T\pa_wV-$$$$\frac{1}{\mathbf{m}}(\pa_vV\times N_0+\mathbf{m}\pa_vN_0)^T\pa_u(V+x_0)V^T\pa_wV+
\frac{|V|^2}{\mathbf{m}}N_0^T(\pa_uN_0\times\pa_vN_0)V^T\pa_wV]]]+$$$$
\frac{N_0^T[\pa_wV\times d(V+x_0)]}{N_0^T[\pa_u(V+x_0)\times\pa_v(V+x_0)]}\odot\frac{d(V+x_0)^TV}{\mathbf{m}^2}
[N_0^T(\pa_wV\times V)N_0^T[\pa_u(V+x_0)\times\pa_v(V+x_0)]+$$$$|V|^2N_0^T(\pa_wV\times V)N_0^T(\pa_uN_0\times\pa_vN_0)-N_0^T(\pa_wV\times V)[N_0^T[\pa_uV\times\pa_v(V+x_0)]-$$$$N_0^T[\pa_vV\times\pa_u(V+x_0)]]
+(N_0\times\pa_wV)^T\pa_uV(N_0\times V)^T\pa_v(V+x_0)-(N_0\times\pa_wV)^T\pa_vV(N_0\times V)^T\pa_u(V+x_0)+$$$$\mathbf{m}^2N_0^T(\pa_wV\times V)N_0^T(\pa_uN_0\times\pa_vN_0)]+
\frac{N_0^T[\pa_wV\times d(V+x_0)]}{N_0^T[\pa_u(V+x_0)\times\pa_v(V+x_0)]}\odot\frac{d(V+x_0)^T(V\times N_0)}{\mathbf{m}^2}[$$$$|V|^2V^T\pa_wVN_0^T(\pa_uN_0\times\pa_vN_0)-V^T\pa_wV[N_0^T[\pa_uV\times\pa_v(V+x_0)]
-N_0^T[\pa_vV\times\pa_u(V+x_0)]]+$$$$\pa_wV^T\pa_uV(N_0\times V)^T\pa_v(V+x_0)
-\pa_wV^T\pa_vV(N_0\times V)^T\pa_u(V+x_0)-$$$$[(\pa_uV\times V)\times(N_0\times\pa_wV)]^T\pa_v(V+x_0)
+[(\pa_vV\times V)\times(N_0\times\pa_wV)]^T\pa_u(V+x_0)]+$$$$2\frac{d(V+x_0)^T\odot dN_0}{\mathbf{m}}N_0^T(\pa_wV\times V)\pa_wV^TV-\frac{1}{\mathbf{m}}d(V+x_0)^TV\odot[dN_0^T\pa_wVN_0^T(\pa_wV\times V)-$$$$
\frac{1}{\mathbf{m}}(\pa_wV\times N_0)^Td(V+x_0)N_0^T(\pa_wV\times V)]-\frac{1}{\mathbf{m}}d(V+x_0)^T(V\times N_0)\odot[$$$$-\frac{1}{\mathbf{m}}dx_0^TV\frac{N_0^T[\pa_wV\times\pa_u(V+x_0)]}{\pa_ux_0^T(N_0\times V)}
N_0^T(\pa_wV\times V)+dN_0^T\pa_wVV^T\pa_wV-\frac{1}{\mathbf{m}}(\pa_wV\times N_0)^Td(V+x_0)V^T\pa_wV]+$$$$\frac{N_0^T[\pa_wV\times d(V+x_0)]}{N_0^T[\pa_u(V+x_0)\times\pa_v(V+x_0)]}\odot
[\frac{d(V+x_0)^T\pa_uN_0}{\mathbf{m}}[V^T\pa_v(V+x_0)N_0^T(\pa_wV\times V)+$$$$(V\times N_0+\mathbf{m}N_0)^T\pa_v(V+x_0)V^T\pa_wV]-\frac{d(V+x_0)^T\pa_vN_0}{\mathbf{m}}[V^T\pa_u(V+x_0)N_0^T(\pa_wV\times V)+$$$$(V\times N_0+\mathbf{m}N_0)^T\pa_u(V+x_0)V^T\pa_wV]-\frac{1}{\mathbf{m}}d(V+x_0)^TV[
-\frac{1}{\mathbf{m}}N_0^T(\pa_wV\times\pa_uV)(V\times N_0)^T\pa_v(V+x_0)+$$$$\frac{1}{\mathbf{m}}N_0^T(\pa_wV\times\pa_vV)(V\times N_0)^T\pa_u(V+x_0)
+\frac{1}{\mathbf{m}}(\pa_uV\times N_0+\mathbf{m}\pa_uN_0)^T\pa_v(V+x_0)N_0^T(\pa_wV\times V)-$$$$
\frac{1}{\mathbf{m}}(\pa_vV\times N_0+\mathbf{m}\pa_vN_0)^T\pa_u(V+x_0)N_0^T(\pa_wV\times V)
+\frac{|V|^2}{\mathbf{m}}N_0^T(\pa_uN_0\times\pa_vN_0)N_0^T(\pa_wV\times V)]-$$$$
\frac{1}{\mathbf{m}}d(V+x_0)^T(V\times N_0)[(-\frac{1}{\mathbf{m}}\pa_ux_0^TV\frac{N_0^T[\pa_wV\times\pa_u(V+x_0)]}{\pa_ux_0^T(N_0\times V)}
+\frac{1}{\mathbf{m}}\pa_wV^T\pa_u(V+x_0)-$$$$\frac{1}{\mathbf{m}}\pa_uV^T\pa_wV)(V\times N_0)^T\pa_v(V+x_0)
-\frac{1}{\mathbf{m}}N_0^T(\pa_wV\times\pa_uV)V^T\pa_v(V+x_0)+$$$$
(\frac{1}{\mathbf{m}}\pa_vx_0^TV\frac{N_0^T[\pa_wV\times\pa_u(V+x_0)]}{\pa_ux_0^T(N_0\times V)}
-\frac{1}{\mathbf{m}}\pa_wV^T\pa_v(V+x_0)+$$$$\frac{1}{\mathbf{m}}\pa_vV^T\pa_wV)(V\times N_0)^T\pa_u(V+x_0)
+\frac{1}{\mathbf{m}}N_0^T(\pa_wV\times\pa_vV)V^T\pa_u(V+x_0)+
\frac{1}{\mathbf{m}}\pa_uV^T\pa_v(V+x_0)N_0^T(\pa_wV\times V)+$$$$\frac{1}{\mathbf{m}}(\pa_uV\times N_0+\mathbf{m}\pa_uN_0)^T\pa_v(V+x_0)V^T\pa_wV-\frac{1}{\mathbf{m}}\pa_vV^T\pa_u(V+x_0)N_0^T(\pa_wV\times V)-$$$$\frac{1}{\mathbf{m}}(\pa_vV\times N_0+\mathbf{m}\pa_vN_0)^T\pa_u(V+x_0)V^T\pa_wV+
\frac{|V|^2}{\mathbf{m}}N_0^T(\pa_uN_0\times\pa_vN_0)V^T\pa_wV]]\}=0.$$

In the above equation the terms containing the second fundamental form of $x_0$ that appear quadratically become, via Gau\ss's theorem, dependent on the metric of $x_0$ and the terms containing the second fundamental of $x_0$ that appear linearly disappear (note $N_0^Td(V+x_0)=-dN_0^TV$) and we get only:

$$\mathbf{n}\{\mathbf{n}^2[\frac{N_0^T[\pa_wV\times d(V+x_0)]}{N_0^T[\pa_u(V+x_0)\times\pa_v(V+x_0)]}\odot\frac{d(V+x_0)^TV}{\mathbf{m}^2}
[N_0^T(\pa_wV\times V)N_0^T[\pa_u(V+x_0)\times\pa_v(V+x_0)]-$$$$N_0^T(\pa_wV\times V)N_0^T(\pa_uV\times\pa_vV)+(\pa_wV\times\pa_vx_0)^TN_0N_0^T(V\times\pa_uV)-
$$$$(\pa_wV\times\pa_ux_0)^TN_0N_0^T(V\times\pa_vV)]+\frac{1}{\mathbf{m}^2}d(V+x_0)^TV\odot dx_0^TV\frac{N_0^T[\pa_wV\times\pa_u(V+x_0)]}{\pa_ux_0^T(N_0\times V)}V^T\pa_wV+$$$$
\frac{N_0^T[\pa_wV\times d(V+x_0)]}{N_0^T[\pa_u(V+x_0)\times\pa_v(V+x_0)]}\odot\frac{d(V+x_0)^TV}{\mathbf{m}^2}
[[\pa_ux_0^TV\frac{N_0^T[\pa_wV\times\pa_u(V+x_0)]}{\pa_ux_0^T(N_0\times V)}-$$$$\pa_wV^T\pa_ux_0]V^T\pa_v(V+x_0)-
[\pa_vx_0^TV\frac{N_0^T[\pa_wV\times\pa_u(V+x_0)]}{\pa_ux_0^T(N_0\times V)}-\pa_wV^T\pa_vx_0]$$$$V^T\pa_u(V+x_0)]]+\mathbf{n}[\frac{N_0^T[\pa_wV\times d(V+x_0)]}{N_0^T[\pa_u(V+x_0)\times\pa_v(V+x_0)]}\odot\frac{d(V+x_0)^TV}{\mathbf{m}^2}[
|V|^2V^T\pa_wVN_0^T(\pa_uN_0\times\pa_vN_0)$$$$-V^T\pa_wV[N_0^T[\pa_uV\times\pa_v(V+x_0)]
-N_0^T[\pa_vV\times\pa_u(V+x_0)]]+\pa_wV^T\pa_uV(N_0\times V)^T\pa_v(V+x_0)-$$$$\pa_wV^T\pa_vV(N_0\times V)^T\pa_u(V+x_0)-\pa_uV^T(N_0\times\pa_wV)V^T\pa_v(V+x_0)+\pa_vV^T(N_0\times\pa_wV)V^T\pa_u(V+x_0)+$$$$
N_0^T(\pa_wV\times V)(\pa_uV^T\pa_vx_0-\pa_vV^T\pa_ux_0)]+$$$$\frac{N_0^T[\pa_wV\times d(V+x_0)]}{N_0^T[\pa_u(V+x_0)\times\pa_v(V+x_0)]}\odot\frac{d(V+x_0)^T(V\times N_0)}{\mathbf{m}^2}[N_0^T(\pa_wV\times V)N_0^T[\pa_u(V+x_0)\times\pa_v(V+x_0)]-$$$$
\pa_uV^T\pa_wVV^T\pa_v(V+x_0)+\pa_vV^T\pa_wVV^T\pa_u(V+x_0)+V^T\pa_wV(\pa_uV^T\pa_vx_0-\pa_vV^T\pa_ux_0)]
+$$$$\frac{1}{\mathbf{m}^2}d(V+x_0)^TV\odot dx_0^TV\frac{N_0^T[\pa_wV\times\pa_u(V+x_0)]}{\pa_ux_0^T(N_0\times V)}
N_0^T(\pa_wV\times V)+$$$$\frac{1}{\mathbf{m}^2}d(V+x_0)^T(V\times N_0)\odot dx_0^TV\frac{N_0^T[\pa_wV\times\pa_u(V+x_0)]}{\pa_ux_0^T(N_0\times V)}V^T\pa_wV+$$$$\frac{1}{\mathbf{m}^2}d(V+x_0)^TV\odot (\pa_wV\times N_0)^Td(V+x_0)V^T\pa_wV+\frac{N_0^T[\pa_wV\times d(V+x_0)]}{N_0^T[\pa_u(V+x_0)\times\pa_v(V+x_0)]}\odot[
-$$$$
\frac{1}{\mathbf{m}}d(V+x_0)^T(V\times N_0)[(-\frac{1}{\mathbf{m}}\pa_ux_0^TV\frac{N_0^T[\pa_wV\times\pa_u(V+x_0)]}{\pa_ux_0^T(N_0\times V)}+\frac{1}{\mathbf{m}}\pa_wV^T\pa_ux_0)V^T\pa_v(V+x_0)$$$$
+(\frac{1}{\mathbf{m}}\pa_vx_0^TV\frac{N_0^T[\pa_wV\times\pa_u(V+x_0)]}{\pa_ux_0^T(N_0\times V)}-\frac{1}{\mathbf{m}}\pa_wV^T\pa_vx_0)V^T\pa_u(V+x_0)$$$$
+\frac{1}{\mathbf{m}}V^T\pa_wV(\pa_uV^T\pa_vx_0-\pa_vV^T\pa_ux_0)]-$$$$
\frac{1}{\mathbf{m}}d(V+x_0)^TV[(-\frac{1}{\mathbf{m}}\pa_ux_0^TV\frac{N_0^T[\pa_wV\times
\pa_u(V+x_0)]}{\pa_ux_0^T(N_0\times V)}+\frac{1}{\mathbf{m}}\pa_wV^T\pa_ux_0)(V\times N_0)^T\pa_v(V+x_0)+$$$$(\frac{1}{\mathbf{m}}\pa_vx_0^TV\frac{N_0^T[\pa_wV\times\pa_u(V+x_0)]}{\pa_ux_0^T(N_0\times V)}-\frac{1}{\mathbf{m}}\pa_wV^T\pa_vx_0)(V\times N_0)^T\pa_u(V+x_0)-$$$$\frac{1}{\mathbf{m}}N_0^T(\pa_wV\times\pa_uV)V^T\pa_v(V+x_0)+
\frac{1}{\mathbf{m}}N_0^T(\pa_wV\times\pa_vV)V^T\pa_u(V+x_0)+$$$$
\frac{1}{\mathbf{m}}N_0^T(\pa_wV\times V)(\pa_uV^T\pa_vx_0-\pa_vV^T\pa_ux_0)
+\frac{1}{\mathbf{m}}(\pa_uV\times N_0)^T\pa_v(V+x_0)V^T\pa_wV-$$$$\frac{1}{\mathbf{m}}(\pa_vV\times N_0)^T\pa_u(V+x_0)V^T\pa_wV+
\frac{|V|^2}{\mathbf{m}}N_0^T(\pa_uN_0\times\pa_vN_0)V^T\pa_wV]]]+$$$$
\frac{N_0^T[\pa_wV\times d(V+x_0)]}{N_0^T[\pa_u(V+x_0)\times\pa_v(V+x_0)]}\odot\frac{d(V+x_0)^TV}{\mathbf{m}^2}
[N_0^T(\pa_wV\times V)N_0^T[\pa_u(V+x_0)\times\pa_v(V+x_0)]+$$$$|V|^2N_0^T(\pa_wV\times V)N_0^T(\pa_uN_0\times\pa_vN_0)-N_0^T(\pa_wV\times V)[N_0^T[\pa_uV\times\pa_v(V+x_0)]-$$$$N_0^T[\pa_vV\times\pa_u(V+x_0)]]
+(N_0\times\pa_wV)^T\pa_uV(N_0\times V)^T\pa_v(V+x_0)-(N_0\times\pa_wV)^T\pa_vV(N_0\times V)^T\pa_u(V+x_0)+$$$$\mathbf{m}^2N_0^T(\pa_wV\times V)N_0^T(\pa_uN_0\times\pa_vN_0)]+
\frac{N_0^T[\pa_wV\times d(V+x_0)]}{N_0^T[\pa_u(V+x_0)\times\pa_v(V+x_0)]}\odot\frac{d(V+x_0)^T(V\times N_0)}{\mathbf{m}^2}[$$$$|V|^2V^T\pa_wVN_0^T(\pa_uN_0\times\pa_vN_0)-V^T\pa_wV[N_0^T[\pa_uV\times\pa_v(V+x_0)]
-N_0^T[\pa_vV\times\pa_u(V+x_0)]]+$$$$\pa_wV^T\pa_uV(N_0\times V)^T\pa_v(V+x_0)
-\pa_wV^T\pa_vV(N_0\times V)^T\pa_u(V+x_0)-$$$$
\pa_uV^T(N_0\times\pa_wV)V^T\pa_v(V+x_0)+\pa_vV^T(N_0\times\pa_wV)V^T\pa_u(V+x_0)
+N_0^T(\pa_wV\times V)(\pa_uV^T\pa_vx_0-\pa_vV^T\pa_ux_0)]
+$$$$\frac{1}{\mathbf{m}^2}d(V+x_0)^TV\odot(\pa_wV\times N_0)^Td(V+x_0)N_0^T(\pa_wV\times V)-\frac{1}{\mathbf{m}}d(V+x_0)^T(V\times N_0)\odot[$$$$-\frac{1}{\mathbf{m}}dx_0^TV\frac{N_0^T[\pa_wV\times\pa_u(V+x_0)]}{\pa_ux_0^T(N_0\times V)}
N_0^T(\pa_wV\times V)-\frac{1}{\mathbf{m}}(\pa_wV\times N_0)^Td(V+x_0)V^T\pa_wV]+$$$$\frac{N_0^T[\pa_wV\times d(V+x_0)]}{N_0^T[\pa_u(V+x_0)\times\pa_v(V+x_0)]}\odot
[d(V+x_0)^T(N_0\times V)V^T\pa_wVN_0^T(\pa_uN_0\times\pa_vN_0)-$$$$\frac{1}{\mathbf{m}}d(V+x_0)^TV[
-\frac{1}{\mathbf{m}}N_0^T(\pa_wV\times\pa_uV)(V\times N_0)^T\pa_v(V+x_0)+$$$$\frac{1}{\mathbf{m}}N_0^T(\pa_wV\times\pa_vV)(V\times N_0)^T\pa_u(V+x_0)
+\frac{1}{\mathbf{m}}(\pa_uV\times N_0)^T\pa_v(V+x_0)N_0^T(\pa_wV\times V)-$$$$
\frac{1}{\mathbf{m}}(\pa_vV\times N_0)^T\pa_u(V+x_0)N_0^T(\pa_wV\times V)
+\frac{|V|^2}{\mathbf{m}}N_0^T(\pa_uN_0\times\pa_vN_0)N_0^T(\pa_wV\times V)]-$$$$
\frac{1}{\mathbf{m}}d(V+x_0)^T(V\times N_0)[(-\frac{1}{\mathbf{m}}\pa_ux_0^TV\frac{N_0^T[\pa_wV\times\pa_u(V+x_0)]}{\pa_ux_0^T(N_0\times V)}
+\frac{1}{\mathbf{m}}\pa_wV^T\pa_u(V+x_0)-$$$$\frac{1}{\mathbf{m}}\pa_uV^T\pa_wV)(V\times N_0)^T\pa_v(V+x_0)
-\frac{1}{\mathbf{m}}N_0^T(\pa_wV\times\pa_uV)V^T\pa_v(V+x_0)+$$$$
(\frac{1}{\mathbf{m}}\pa_vx_0^TV\frac{N_0^T[\pa_wV\times\pa_u(V+x_0)]}{\pa_ux_0^T(N_0\times V)}
-\frac{1}{\mathbf{m}}\pa_wV^T\pa_v(V+x_0)+$$$$\frac{1}{\mathbf{m}}\pa_vV^T\pa_wV)(V\times N_0)^T\pa_u(V+x_0)
+\frac{1}{\mathbf{m}}N_0^T(\pa_wV\times\pa_vV)V^T\pa_u(V+x_0)+
\frac{1}{\mathbf{m}}\pa_uV^T\pa_vx_0N_0^T(\pa_wV\times V)+$$$$\frac{1}{\mathbf{m}}(\pa_uV\times N_0)^T\pa_v(V+x_0)V^T\pa_wV-\frac{1}{\mathbf{m}}\pa_vV^T\pa_ux_0N_0^T(\pa_wV\times V)-$$$$\frac{1}{\mathbf{m}}(\pa_vV\times N_0)^T\pa_u(V+x_0)V^T\pa_wV+
\frac{|V|^2}{\mathbf{m}}N_0^T(\pa_uN_0\times\pa_vN_0)V^T\pa_wV]]\}=0.$$

\subsection{The case $m=V+\mathbf{m}N_0$}\noindent

\noindent

We have $0=(Rm)^T\ti d(RV+x)=m^T(\om\times V+d(V+x_0)+\pa_wVdw)$, or, assuming $m^T\pa_wV\neq 0$,

$$dw=\frac{\mathbf{m}V^T(\om\times N_0+dN_0)-V^Td(V+x_0)}
{V^T\pa_wV}.$$

Applying the compatibility condition $d$ to this equation and using the equation itself we get

$$0=-\frac{\pa_w[\mathbf{m}V^T(\om\times N_0+dN_0)-V^Td(V+x_0)]}
{V^T\pa_wV}\wedge\frac{\mathbf{m}V^T(\om\times N_0+dN_0)-V^Td(V+x_0)}
{V^T\pa_wV}+$$
$$d\frac{\mathbf{m}V^T}{V^T\pa_wV}\wedge(\om\times N_0+dN_0)-d\frac{V^T}{V^T\pa_wV}\wedge d(V+x_0)
\stackrel{(\ref{eq:fund})}{=}$$
$$-\frac{\pa_w[V^Td(V+x_0)]}{V^T\pa_wV}\wedge\frac{V^Td(V+x_0)}{V^T\pa_wV}-d\frac{V^T}{V^T\pa_wV}\wedge d(V+x_0)-\frac{\mathbf{m}^2N_0^T(\pa_wV\times V)N_0^T(dN_0\times\wedge dN_0)}{2(V^T\pa_wV)^2}$$
$$+[-\frac{\pa_w\mathbf{m}V^Td(V+x_0)}{(V^T\pa_wV)^2}+\frac{\mathbf{m}\pa_w[V^Td(V+x_0)]}{(V^T\pa_wV)^2}
+d\frac{\mathbf{m}}{V^T\pa_wV}]\wedge V^T(\om\times N_0+dN_0)$$
$$+\frac{\mathbf{m}N_0^T[\pa_wV\times d(V+x_0)]}{(V^T\pa_wV)^2}\wedge(N_0\times V)^T(\om\times N_0+dN_0).$$

As in the previous case we get

$$\frac{\mathbf{m}N_0^T[\pa_wV\times d(V+x_0)]}{(V^T\pa_wV)^2}=\frac{bdx_0^T(N_0\times V)}{|N_0\times V|^2}$$
with $b\subset\mathbb{C}$.


\begin{thebibliography}{99}
\def\topset{0pt}
\def\parsep{0pt plus 5pt minus 1pt}
\def\itemsep{-0.5ex}
\small

\bibitem{B3} L. Bianchi {\it Concerning Singular Transformations $B_k$ of surfaces applicable
to quadrics,} Transactions of the American Mathematical Society,
{\href{http://www.ams.org/journals/tran/1917-018-03/S0002-9947-1917-1501075-2/S0002-9947-1917-1501075-2.pdf}
{{\bf 18} (1917), 379-401}}.

\bibitem{D1} G. Darboux  {\href{http://fr.dleex.com/details/?11194}
{\it Le\c{c}ons Sur La Th\'{e}orie G\'{e}n\'{e}rale Des
Surfaces,}} Vol {\bf 1-4}, Gauthier-Villars, Paris (1894-1917).

\bibitem{D3} I. Dinc\u{a} {\it On Bianchi's B\"{a}cklund transformation of quadrics,}
{\href{http://arxiv.org/abs/1110.5474} {arxiv:1110.5474v2}} and Journal of Geometry and Physics, {\bf 73} (2013), 104-124.

\end{thebibliography}
\end{document}